# FRACTIONAL CAUCHY PROBLEMS ON BOUNDED DOMAINS[1]

By Mark M. Meerschaert[2], Erkan Nane and P. Vellaisamy

*Michigan State University, Auburn University and Indian Institute of Technology Bombay*

Fractional Cauchy problems replace the usual first-order time derivative by a fractional derivative. This paper develops classical solutions and stochastic analogues for fractional Cauchy problems in a bounded domain $D \subset \mathbb{R}^d$ with Dirichlet boundary conditions. Stochastic solutions are constructed via an inverse stable subordinator whose scaling index corresponds to the order of the fractional time derivative. Dirichlet problems corresponding to iterated Brownian motion in a bounded domain are then solved by establishing a correspondence with the case of a half-derivative in time.

**1. Introduction.** In this paper, we extend the approach of Meerschaert and Scheffler [23] and Meerschaert et al. [24] to fractional Cauchy problems on bounded domains. Our methods involve eigenfunction expansions, killed Markov processes and inverse stable subordinators. In a recent related paper [7], we establish a connection between fractional Cauchy problems with index $\beta = 1/2$ on an unbounded domain, and iterated Brownian motion (IBM), defined as $Z_t = B(|Y_t|)$, where $B$ is a Brownian motion with values in $\mathbb{R}^d$ and $Y$ is an independent one-dimensional Brownian motion. Since IBM is also the stochastic solution to a Cauchy problem involving a fourth-order derivative in space [2, 14], that paper also establishes a connection between certain higher-order Cauchy problems and their time-fractional analogues. More generally, Baeumer, Meerschaert and Nane [7] shows a connection between fractional Cauchy problems with $\beta = 1/2$ and higher-order Cauchy problems that involve the square of the generator. In the present paper, we

Received February 2008; revised July 2008.
[1]This paper was completed while E. Nane and P. Vellaisamy were visiting Michigan State University.
[2]Supported in part by NSF Grant DMS-07-06440.

clarify that relationship by identifying the boundary conditions that make these two formulations have the same unique solutions on bounded domains; see Corollary 4.5.

A celebrated paper of Einstein [15] established a mathematical link between random walks, the diffusion equation and Brownian motion. The scaling limits of a simple random walk with mean zero, finite variance jumps yields a Brownian motion. The probability densities of the Brownian motion variables solve a diffusion equation, and hence we refer to the Brownian motion as the stochastic solution to the diffusion equation. The interplay between stochastic processes and partial differential equations has many practical and theoretical applications. For example, it has lead to efficient particle tracking simulation methods for solving complicated boundary value problems [18]. The diffusion equation is the most familiar Cauchy problem. The general abstract Cauchy problem is $\partial_t u = Lu$, where $u(t)$ takes values in a Banach space and $L$ is the generator of a continuous semigroup on that space [4]. If $L$ generates a Markov process, then we call this Markov process a stochastic solution to the Cauchy problem $\partial_t u = Lu$, since its probability densities (or distributions) solve the Cauchy problem. This point of view has proven useful, for instance, in the modern theory of fractional calculus, since fractional derivatives are generators of certain ($\alpha$-stable) stochastic processes [23].

Fractional derivatives are almost as old as their more familiar integer-order counterparts [26, 30]. Fractional diffusion equations have recently been applied to problems in physics, finance, hydrology and many other areas [17, 20, 25, 32]. Fractional space derivatives are used to model anomalous diffusion or dispersion, where a particle plume spreads at a rate inconsistent with the classical Brownian motion model, and the plume may be asymmetric. When a fractional derivative replaces the second derivative in a diffusion or dispersion model, it leads to enhanced diffusion (also called superdiffusion). Fractional time derivatives are connected with anomalous subdiffusion, where a cloud of particles spreads more slowly than a classical diffusion. Fractional Cauchy problems replace the integer time derivative by its fractional counterpart: $\partial_t^\beta u = Lu$. Here, $\partial_t^\beta g(t)$ indicates the Caputo fractional derivative in time, the inverse Laplace transform of $s^\beta \tilde{g}(s) - s^{\beta-1} g(0)$, where $\tilde{g}(s) = \int_0^\infty e^{-st} g(t)\,dt$ is the usual Laplace transform [11]. Nigmatullin [27] gave a physical derivation of the fractional Cauchy problem, when $L$ is the generator of some continuous Markov process $\{Y(t)\}$ started at $x = 0$. The mathematical study of fractional Cauchy problems was initiated by [19, 20, 33]. The existence and uniqueness of solutions was proved in [19, 20]. Fractional Cauchy problems were also invented independently by Zaslavsky [35] as a model for Hamiltonian chaos.

Stochastic solutions of fractional Cauchy problems are subordinated processes. If $X(t)$ is a stochastic solution to the Cauchy problem $\partial_t u = Au$, then under certain technical conditions, the subordinate process $X(E_t)$ is a stochastic solution to the fractional Cauchy problem $\partial_t^\beta u = Au$; see [5]. Here, $E_t$ is the inverse or hitting time process to a stable subordinator $D_t$ with index $\beta \in (0,1)$. That is, $E_t = \inf\{x > 0 \colon D_x > t\}$, and $D_t$ is a Lévy process (continuous in probability with independent, stationary increments) whose smooth probability density $g_\beta(t)$ has Laplace transform $e^{-s^\beta} = \tilde{g}_\beta(s)$; see [31]. Just as Brownian motion is a scaling limit of a simple random walk, the stochastic solution to certain fractional Cauchy problems are scaling limits of continuous time random walks, in which the i.i.d. jumps are separated by i.i.d. waiting times [24]. Fractional time derivatives arise from power law waiting times, where the probability of waiting longer than time $t > 0$ falls off like $t^{-\beta}$ for $t$ large [23]. This is related to the fact that fractional derivatives are nonlocal operators defined by convolution with a power law [5].

**2. Preliminaries.** Let $D$ be a bounded domain in $\mathbb{R}^d$. We define the following spaces of functions:

$$C(D) = \{u \colon D \to \mathbb{R} \colon u \text{ is continuous}\},$$
$$C(\bar{D}) = \{u \colon \bar{D} \to \mathbb{R} \colon u \text{ is uniformly continuous}\},$$
$$C^k(D) = \{u \colon D \to \mathbb{R} \colon u \text{ is } k\text{-times continuously differentiable}\},$$
$$C^k(\bar{D}) = \{u \in C^k(D) \colon D^\gamma u \text{ is uniformly continuous for all } |\gamma| \le k\}.$$

Thus, if $u \in C^k(\bar{D})$, then $D^\gamma u$ continuously extends to $\bar{D}$ for each multi-index $\gamma$ with $|\gamma| \le k$.

We define the spaces of functions $C^\infty(D) = \bigcap_{k=1}^\infty C^k(D)$ and $C^\infty(\bar{D}) = \bigcap_{k=1}^\infty C^k(\bar{D})$.

Also, let $C^{k,\alpha}(D)$ $[C^{k,\alpha}(\bar{D})]$ be the subspace of $C^k(D)$ $[C^k(\bar{D})]$ that consists of functions whose $k$th order partial derivatives are uniformly Hölder continuous with exponent $\alpha$ in $D$. For simplicity, we will write

$$C^{0,\alpha}(D) = C^\alpha(D), \qquad C^{0,\alpha}(\bar{D}) = C^\alpha(\bar{D})$$

with the understanding that $0 < \alpha < 1$ whenever this notation is used, unless otherwise stated.

We use $C_c(D), C_c^k(D), C_c^{k,\alpha}(D)$ to denote those functions in $C(D), C^k(D), C^{k,\alpha}(D)$ with compact support.

A subset $D$ of $\mathbb{R}^d$ is an $l$-dimensional manifold with boundary if every point of $D$ possesses a neighborhood diffeomorphic to an open set in the space $H^l$, which is the upper half space in $\mathbb{R}^l$. Such a diffeomorphism is called a local parametrization of $D$. The boundary of $D$, denoted by $\partial D$, consists of those points that belong to the image of the boundary of $H^l$ under some local parametrization. If the diffeomorphism and its inverse are $C^{k,\alpha}$ functions, then we write $\partial D \in C^{k,\alpha}$.

We say that $D$ satisfies an *exterior cone condition* at a fixed point $x_0 \in \partial D$ if there exists a finite right circular cone $V = V_{x_0}$ with vertex $x_0$ such that $\bar{D} \cap V_{x_0} = x_0$, and a *uniform exterior cone condition* if $D$ satisfies an exterior cone condition at every point $x_0 \in \partial D$ and the cones $V_{x_0}$ are all congruent to some fixed cone $V$.

Let $\{X_t\}$ be a Brownian motion in $\mathbb{R}^d$ and $A \subset \mathbb{R}^d$. Let $T_A = \inf\{t > 0 : X_t \in A\}$ be the first hitting time of the set $A$. We say that a point $y$ is *regular* for a set $A$ if $P_y[T_A = 0] = 1$. Note that a point $y$ is regular for a set $A$ provided, starting at $y$, the process does not go a positive length of time before hitting $A$.

The right condition for the existence of the solution to the Dirichlet problem turns out to be that every point of $\partial D$ is regular for $D^C$ (cf. [8], Section II.1).

If a domain satisfies a uniform exterior cone condition, then every point of $\partial D$ is regular for $D^C$.

The Sobolev space $W^{k,p}(D)$ consists of all locally integrable functions $u : D \to \mathbb{R}$ such that for each multi-index $\gamma$ with $|\gamma| \leq k$, $D^\gamma u$ exists in the weak sense and belongs to $L^p(D)$ with the following norm:

$$\|f\|_{W^{k,p}(D)} = \left( \|f\|_{p,D}^p + \sum_{|\gamma| \leq k} \|D^\gamma f\|_{p,D}^p \right)^{1/p}.$$

We denote by $W_0^{k,p}(D)$ the closure of $C_c^k(D)$ in $W^{k,p}(D)$. If $p = 2$, we usually write

$$H^k(D) = W^{k,2}(D), \qquad H_0^k(D) = W_0^{k,2}(D) \qquad (k = 0, 1, \ldots).$$

Let $\partial D \in C^1$. Then at each point $x \in \partial D$ there exists a unique outward pointing unit vector

$$\theta(x) = (\theta_1(x), \ldots, \theta_d(x)).$$

Let $u \in C^1(\bar{D})$, the set of functions which have continuous extension of the first derivative up to the boundary. Let

$$D_\theta u = \frac{\partial u}{\partial \theta} = \theta \cdot \nabla u,$$

denote the directional derivative, where $\nabla u$ is the gradient vector of $u$.

Now we recall Green's first and second identities (see, e.g., [16], Section 2.4). Let $u, v \in C^2(D) \cap C^1(\bar{D})$. Then

$$\int_D \frac{\partial u}{\partial x_i} v \, dx = \int_{\partial D} uv\theta_i \, ds - \int_D u \frac{\partial v}{\partial x_i} \, dx$$
(integration by parts formula),

$$\int_D \nabla v \cdot \nabla u \, dx = -\int_D u \Delta v \, dx + \int_{\partial D} \frac{\partial v}{\partial \theta} u \, ds \quad \text{(Green's first identity)},$$

$$\int_D [u \Delta v - v \Delta u] \, dx = \int_{\partial D} \left[ u \frac{\partial v}{\partial \theta} - v \frac{\partial u}{\partial \theta} \right] ds \quad \text{(Green's second identity)}.$$

A uniformly elliptic operator of divergence form is defined on $C^2$ functions by

$$(2.1) \qquad Lu = \sum_{i,j=1}^{d} \frac{\partial(a_{ij}(x)(\partial u/\partial x_i))}{\partial x_j}$$

with $a_{ij}(x) = a_{ji}(x)$ and, for some $\lambda > 0$,

$$(2.2) \qquad \lambda \sum_{i=1}^{n} y_i^2 \leq \sum_{i,j=1}^{n} a_{ij}(x) y_i y_j \leq \lambda^{-1} \sum_{i=1}^{n} y_i^2 \qquad \forall y \in \mathbb{R}^d.$$

The operator $L$ acts on the Hilbert space $L^2(D)$. We define the initial domain $C_0^\infty(\bar{D})$ of the operator as follows. We say that $f$ is in $C_0^\infty(\bar{D})$, if $f \in C^\infty(\bar{D})$ and $f(x) = 0$ for all $x \in \partial D$. This condition incorporates the notion of Dirichlet boundary conditions.

From [13], Corollary 6.1, we have that the associated quadratic form

$$Q(f,g) = \int_D \sum_{i,j=1}^{d} a_{ij} \frac{\partial f}{\partial x_i} \frac{\partial g}{\partial x_j} \, dx$$

is closable on the domain $C_0^\infty(\bar{D})$ and the domain of the closure is independent of the particular coefficients $(a_{ij})$ chosen. In particular, for $f, g \in C_0^\infty(\bar{D})$ by integration by parts

$$\int_D g(x) L f(x) \, dx = Q(f,g) = \int_D f(x) L g(x) \, dx,$$

which shows that $L$ is symmetric.

From now on, we will use the symbol $L_D$ if we particularly want to emphasize the choice of Dirichlet boundary conditions, to refer to the self-adjoint operator associated with the closure of the quadratic form above by the use of [13], Theorem 4.4.5. Thus, $L_D$ is the Friedrichs extension of the operator defined initially on $C_0^\infty(\bar{D})$.

If the coefficients $a_{ij}(x)$ are smooth $[a_{ij}(x) \in C^1(D)]$, then $L_D u$ takes the form

$$L_D u = \sum_{i,j=1}^{d} a_{ij}(x) \frac{\partial^2 u}{\partial x_i \, \partial x_j} + \sum_{i=1}^{d} \left( \sum_{j=1}^{d} \frac{\partial a_{ij}(x)}{\partial x_j} \right) \frac{\partial u}{\partial x_i}$$

$$= \sum_{i,j=1}^{d} a_{ij}(x) \frac{\partial^2 u}{\partial x_i \, \partial x_j} + \sum_{i=1}^{d} b_i(x) \frac{\partial u}{\partial x_i}.$$

If $X_t$ is a solution to

$$dX_t = \sigma(X_t) \, dW_t + b(X_t) \, dt, \qquad X_0 = x_0,$$

where $\sigma$ is a $d \times d$ matrix, and $W_t$ is a Brownian motion, then $X_t$ is associated with the operator $L_D$ with $a = \sigma \sigma^T$ (see Chapters 1 and 5 of Bass [8]). Define the first exit time as $\tau_D(X) = \inf\{t \geq 0 : X_t \notin D\}$. The semigroup defined by $T(t)f(x) = E_x[f(X_t) I(\tau_D(X)) > t)]$ has generator $L_D$, which follows by an application of the Itô formula.

Let $D$ be a bounded domain in $\mathbb{R}^d$. Suppose $L$ is a uniformly elliptic operator of divergence form with Dirichlet boundary conditions on $D$, and that there exists a constant $\Lambda$ such that for all $x \in D$,

$$(2.3) \qquad \sum_{i,j=1}^{d} |a_{ij}(x)| \leq \Lambda.$$

Let $T_D(t)$ be the corresponding semigroup. Then $T_D(t)$ is an ultracontractive semigroup (even an intrinsically ultracontractive semigroup); see Corollary 3.2.8, Theorems 2.1.4 and 4.2.4 and Note 4.6.10 in [12]. Every ultracontractive semigroup has a kernel for the killed semigroup on a bounded domain which can be represented as a series expansion of the eigenvalues and the eigenfunctions of $L_D$ (cf. [12], Theorems 2.1.4 and 2.3.6 and [16], Theorems 8.37 and 8.38): There exist eigenvalues $0 < \mu_1 < \mu_2 \leq \mu_3 \leq \cdots$, such that $\mu_n \to \infty$, as $n \to \infty$, with the corresponding complete orthonormal set (in $H_0^2$) of eigenfunctions $\psi_n$ of the operator $L_D$ satisfying

$$(2.4) \qquad L_D \psi_n(x) = -\mu_n \psi_n(x), \qquad x \in D : \psi_n|_{\partial D} = 0.$$

In this case,

$$p_D(t, x, y) = \sum_{n=1}^{\infty} e^{-\mu_n t} \psi_n(x) \psi_n(y)$$

is the heat kernel of the killed semigroup $T_D$. The series converges absolutely and uniformly on $[t_0, \infty) \times D \times D$ for all $t_0 > 0$.

Denote the Laplace transform by

$$\tilde{u}(s,x) = \int_0^\infty e^{-st} u(t,x)\, dt.$$

Since we are working on a bounded domain, the Fourier transform methods in [24] are not useful. Instead we will employ Hilbert space methods. Hence, given a complete orthonormal basis $\{\psi_n(x)\}$ on $H_0^2$, we will call

$$\bar{u}(t,n) = \int_D \psi_n(x) u(t,x)\, dx,$$

$$\hat{u}(s,n) = \int_D \psi_n(x) \int_0^\infty e^{-st} u(t,x)\, dt\, dx$$

$$= \int_D \psi_n(x) \tilde{u}(s,x)\, dx$$

the $\psi_n$ and $\psi_n$-Laplace transforms, respectively. Since $\{\psi_n\}$ is a complete orthonormal basis for $H_0^2$, we can invert the $\psi_n$-transform

$$u(t,x) = \sum_n \bar{u}(t,n) \psi_n(x)$$

for any $t > 0$, where the sum converges in the $L^2$ sense (e.g., see [29], Proposition 10.8.27).

Suppose $D$ satisfies a uniform exterior cone condition. Let $\{X_t\}$ be a Markov process in $\mathbb{R}^d$ with generator $L_D$, and $f$ be continuous on $\bar{D}$. Then the semigroup

(2.5)
$$\begin{aligned} T_D(t) f(x) &= E_x[f(X_t) I(t < \tau_D(X))] \\ &= \int_D p_D(t,x,y) f(y)\, dy \\ &= \sum_{n=1}^\infty e^{-\mu_n t} \psi_n(x) \bar{f}(n) \end{aligned}$$

solves the Dirichlet initial-boundary value problem in $D$:

$$\frac{\partial u(t,x)}{\partial t} = L_D u(t,x), \qquad x \in D, t > 0,$$
$$u(t,x) = 0, \qquad x \in \partial D,$$
$$u(0,x) = f(x), \qquad x \in D.$$

REMARK 2.1. The eigenfunctions belong to $L^\infty(D) \cap C^\alpha(D)$ for some $\alpha > 0$, by [16], Theorems 8.15 and 8.24. If $D$ satisfies a *uniform exterior cone condition* all the eigenfunctions belong to $C^\alpha(\bar{D})$ by [16], Theorem 8.29. If $a_{ij} \in C^\alpha(\bar{D})$ and $\partial D \in C^{1,\alpha}$, then all the eigenfunctions belong to $C^{1,\alpha}(\bar{D})$ by [16], Corollary 8.36. If $a_{ij} \in C^\infty(D)$ then each eigenfunction of $L$ is in $C^\infty(D)$ by [16], Corollary 8.11. If $a_{ij} \in C^\infty(\bar{D})$ and $\partial D \in C^\infty$, then each eigenfunction of $L$ is in $C^\infty(\bar{D})$ by [16], Theorem 8.13.

([16], Corollary 8.36). Let $D$ be bounded and every point of $\partial D$ be regular for $D^C$. In the case $L_D = \Delta_D$, the corresponding Markov process is a killed Brownian motion. We denote the eigenvalues and the eigenfunctions of $\Delta_D$ by $\{\lambda_n, \phi_n\}_{n=1}^\infty$, where $\phi_n \in C^\infty(D)$. The corresponding heat kernel is given by

$$p_D(t,x,y) = \sum_{n=1}^\infty e^{-\lambda_n t} \phi_n(x) \phi_n(y).$$

The series converges absolutely and uniformly on $[t_0, \infty) \times D \times D$ for all $t_0 > 0$. In this case, the semigroup given by

$$T_D(t)f(x) = E_x[f(X_t) I(t < \tau_D(X))]$$
$$= \int_D p_D(t,x,y) f(y)\, dy = \sum_{n=1}^\infty e^{-\lambda_n t} \phi_n(x) \bar{f}(n)$$

solves the heat equation in $D$ with Dirichlet boundary conditions:

$$\frac{\partial u(t,x)}{\partial t} = \Delta u(t,x), \qquad x \in D, t > 0,$$
$$u(t,x) = 0, \qquad x \in \partial D,$$
$$u(0,x) = f(x), \qquad x \in D.$$

The Caputo fractional derivative of order $0 < \beta < 1$, defined by

$$(2.6) \qquad \frac{\partial^\beta u(t,x)}{\partial t^\beta} = \frac{1}{\Gamma(1-\beta)} \int_0^t \frac{\partial u(s,x)}{\partial s} \frac{ds}{(t-s)^\beta},$$

was invented to properly handle initial values [11]. Its Laplace transform $s^\beta \tilde{u}(s,x) - s^{\beta-1} u(0,x)$ incorporates the initial value in the same way as the first derivative. The Caputo derivative has been widely used to solve ordinary differential equations that involve a fractional time derivative [17, 28].

Inverse stable subordinators arise in [23, 24] as scaling limits of continuous time random walks. Let $S(n) = Y_1 + \cdots + Y_n$ a sum of i.i.d. random variables with $EY_n = 0$ and $EY_n^2 < \infty$. The scaling limit $c^{-1/2} S([ct]) \Rightarrow B(t)$ as $c \to \infty$ is a Brownian motion, normal with mean zero and variance proportional to $t$. Consider $Y_n$ to be the random jumps of a particle. If we impose a random waiting time $T_n$ before the $n$th jump $Y_n$, then the position of the particle at time $T_n = J_1 + \cdots + J_n$ is given by $S(n)$. The number of jumps by time $t > 0$ is $N_t = \max\{n : T_n \leq t\}$, so the position of the particle at time $t > 0$ is $S(N_t)$, a subordinated process. If $P(J_n > t) = t^{-\beta} L(t)$ for some $0 < \beta < 1$, where $L(t)$ is slowly varying, then the scaling limit $c^{-1/\beta}T_{[ct]} \Rightarrow D_t$ is a strictly increasing stable Lévy motion with index $\beta$, sometimes called a stable subordinator. The jump times $T_n$ and the number of jumps $N_t$ are inverses $\{N_t \geq x\} = \{T(\lceil x \rceil) \leq t\}$, where $\lceil x \rceil$ is the smallest integer greater than or equal to $x$. It follows that the scaling limits are also inverses $c^{-\beta}N_{ct} \Rightarrow E_t$, where $E_t = \inf\{x : D_x > t\}$, so that $\{E_t \leq x\} = \{D_x \geq t\}$. Since $N_{ct} \approx c^\beta E_t$, the scaling limit of the particle locations is $c^{-\beta/2}S(N_{[ct]}) \approx (c^\beta)^{-1/2}S(c^\beta E_t) \approx B(E_t)$, a Brownian motion subordinated to the inverse or hitting time (or first passage time) process of the stable subordinator $D_t$. The random variable $D_t$ has a smooth density. For properly scaled waiting times, the density of $D_t$ has Laplace transform $e^{-ts^\beta}$ for any $t > 0$, and $D_t$ is identically distributed with $t^{1/\beta}D_1$. Writing $g_\beta(u)$ for the density of $D_1$, it follows that $D_t$ has density $t^{-1/\beta}g_\beta(t^{-1/\beta}u)$ for any $t > 0$. Using the inverse relation $P(E_t \leq x) = P(D_x \geq t)$ and taking derivatives, it follows that $E_t$ has density

$$f_t(x) = t\beta^{-1}x^{-1-1/\beta}g_\beta(tx^{-1/\beta}), \tag{2.7}$$

whose $t \mapsto s$ Laplace transform $s^{\beta-1}e^{-xs^\beta}$ can also be derived from the equation

$$f_t(x) = \frac{d}{dx}P(D_x \geq t) = \frac{d}{dx}\int_t^\infty x^{-1/\beta}g_\beta(x^{-1/\beta}u)\,du$$

by taking Laplace transforms on both sides.

**3. Fractional diffusion in bounded domains.** Fractional Cauchy problems replace the usual first-order time derivative with its fractional analogue. In this section, we prove classical (strong) solutions to fractional Cauchy problems on bounded domains in $\mathbb{R}^d$. We also give an explicit solution formula, based on the solution of the corresponding Cauchy problem. Our methods are inspired by the approach in [24], where Laplace transforms are used to handle the fractional time derivative, and spatial derivative operators (or more generally, pseudo-differential operators) are treated with Fourier transforms. In the present paper, we use an eigenfunction expansion in place of Fourier transforms, since we are operating on a bounded domain. Our first result, Theorem 3.1, is focused on a fractional diffusion $\partial_t^\beta u = \Delta u$, and we lay out all the details of the argument in the most familiar setting. Then, in Theorem 3.6, we use a separation of variables to handle arbitrary generators: $\partial_t^\beta u = Lu$. Along the way, we explicate the stochastic solutions in terms of killed Markov processes.

Let $\beta \in (0,1)$, $D_\infty = (0,\infty) \times D$ and define

$$\mathcal{H}_\Delta(D_\infty) \equiv \left\{ u : D_\infty \to \mathbb{R} : \frac{\partial}{\partial t} u, \frac{\partial^\beta}{\partial t^\beta} u, \Delta u \in C(D_\infty), \right.$$

$$\left. \left| \frac{\partial}{\partial t} u(t,x) \right| \leq g(x) t^{\beta-1}, \ g \in L^\infty(D), \ t > 0 \right\}.$$

We will write $u \in C^k(\bar{D})$ to mean that for each fixed $t > 0$, $u(t, \cdot) \in C^k(\bar{D})$, and $u \in C_b^k(\bar{D}_\infty)$ to mean that $u \in C^k(\bar{D}_\infty)$ and is bounded.

THEOREM 3.1. *Let $0 < \alpha < 1$. Let $D$ be a bounded domain with $\partial D \in C^{1,\alpha}$, and $T_D(t)$ be the killed semigroup of Brownian motion $\{X_t\}$ in $D$. Let $E_t$ be the process inverse to a stable subordinator of index $\beta \in (0,1)$ independent of $\{X_t\}$. Let $f \in D(\Delta_D) \cap C^1(\bar{D}) \cap C^2(D)$ for which the eigenfunction expansion (of $\Delta f$) with respect to the complete orthonormal basis $\{\phi_n : n \in \mathbb{N}\}$ converges uniformly and absolutely. Then the unique (classical) solution of*

(3.1)
$$u \in \mathcal{H}_\Delta(D_\infty) \cap C_b(\bar{D}_\infty) \cap C^1(\bar{D}),$$
$$\frac{\partial^\beta}{\partial t^\beta} u(t,x) = \Delta u(t,x), \qquad x \in D, t > 0,$$
$$u(t,x) = 0, \qquad x \in \partial D, t > 0,$$
$$u(0,x) = f(x), \qquad x \in D$$

*is given by*

$$u(t,x) = E_x[f(X(E_t)) I(\tau_D(X) > E_t)]$$
$$= \frac{t}{\beta} \int_0^\infty T_D(l) f(x) g_\beta(tl^{-1/\beta}) l^{-1/\beta-1} \, dl = \int_0^\infty T_D((t/l)^\beta) f(x) g_\beta(l) \, dl.$$

PROOF. Assume that $u(t,x)$ solves (3.1). Using Green's second identity, we obtain

$$\int_D [u \Delta \phi_n - \phi_n \Delta u] \, dx = \int_{\partial D} \left[ u \frac{\partial \phi_n}{\partial \theta} - \phi_n \frac{\partial u}{\partial \theta} \right] ds = 0,$$

since $u|_{\partial D} = 0 = \phi_n|_{\partial D}$ and $u, \phi_n \in C^1(\bar{D})$. Hence, the $\phi_n$-transform of $\Delta u$ is

$$\int_D \phi_n(x) \Delta u(t,x) \, dx = \int_D u(t,x) \Delta \phi_n(x) \, dx$$
$$= -\lambda_n \int_D u(t,x) \phi_n(x) \, dx = -\lambda_n \bar{u}(t,n),$$

as $\phi_n$ is the eigenfunction of the Laplacian corresponding to eigenvalue $\lambda_n$.

Recall the definition (2.6) of the fractional derivative, and use the fact that $|\frac{\partial u(t,x)}{\partial t}|$ is bounded by $g(x)t^{\beta-1}$ where $g \in L^\infty(D)$. It follows that the fractional derivative commutes with $\phi_n$-transform, and hence, we can take $\phi_n$-transforms in (3.1) to obtain

$$\text{(3.2)} \qquad \frac{\partial^\beta}{\partial t^\beta}\bar{u}(t,n) = -\lambda_n \bar{u}(t,n).$$

Since $u$ is uniformly continuous on $C([0,\epsilon]) \times \bar{D}$ and hence uniformly bounded on $[0,\epsilon] \times \bar{D}$, we have by dominated convergence that $\lim_{t\to 0}\int_D u(t,x) \phi_n(x)\,dx = \bar{f}(n)$. Hence, $\bar{u}(0,n) = \bar{f}(n)$. A similar argument shows that $t \mapsto \bar{u}(t,n)$ is a continuous function of $t \in [0,\infty)$ for every $n$. Then, taking Laplace transforms on both sides of (3.2), we get

$$\text{(3.3)} \qquad s^\beta \hat{u}(s,n) - s^{\beta-1}\bar{u}(0,n) = -\lambda_n \hat{u}(s,n)$$

which leads to

$$\text{(3.4)} \qquad \hat{u}(s,n) = \frac{\bar{f}(n)s^{\beta-1}}{s^\beta + \lambda_n}.$$

By inverting the above Laplace transform, we obtain

$$\bar{u}(t,n) = \bar{f}(n)E_\beta(-\lambda_n t^\beta)$$

in terms of the Mittag–Leffler function defined by

$$E_\beta(z) = \sum_{k=0}^{\infty} \frac{z^k}{\Gamma(1+\beta k)};$$

see, for example, [22]. Inverting the $\phi_n$-transform, we get an $L^2$-convergent solution of equation (3.1) as

$$\text{(3.5)} \qquad u(t,x) = \sum_{n=1}^{\infty} \bar{f}(n)\phi_n(x)E_\beta(-\lambda_n t^\beta)$$

for each $t \geq 0$. In order to complete the proof, it will suffice to show that the series (3.5) converges pointwise, and satisfies all the conditions in (3.1).

STEP 1. We begin showing that (3.5) convergence uniformly in $t \in [0,\infty)$ in the $L^2$ sense. Define the sequence of functions

$$\text{(3.6)} \qquad u_N(t,x) = \sum_{n=1}^{N} \bar{f}(n)\phi_n(x)E_\beta(-\lambda_n t^\beta).$$

We have, by [21], equation (13),

$$\text{(3.7)} \qquad 0 \leq E_\beta(-\lambda_n t^\beta) \leq c/(1+\lambda_n t^\beta).$$

Since $f \in L^2(D)$, we can write $f(x) = \sum_{n=1}^{\infty} \bar{f}(n) \phi_n(x)$, and then the Parseval identity yields

$$\sum_{n=1}^{\infty} (\bar{f}(n))^2 = \|f\|_{2,D}^2 < \infty.$$

Then, given $\epsilon > 0$, we can choose $n_0(\epsilon)$ such that

(3.8) $$\sum_{n=n_0(\epsilon)}^{\infty} (\bar{f}(n))^2 < \epsilon.$$

For $N > M > n_0(\epsilon)$ and $t \geq 0$,

$$\|u_N(t,x) - u_M(t,x)\|_{2,D}^2 \leq \left\| \sum_{n=M}^{N} \bar{f}(n) \phi_n(x) E_\beta(-\lambda_n t^\beta) \right\|_{2,D}^2$$

(3.9) $$\leq \left[ \frac{c}{(1+\lambda_{n_0} t^\beta)} \right]^2 \sum_{n=n_0(\epsilon)} (\bar{f}(n))^2$$

$$\leq c^2 \sum_{n=n_0(\epsilon)} (\bar{f}(n))^2 < c^2 \epsilon.$$

Thus, the series (3.5) converges in $L^2(D)$, uniformly in $t \geq 0$.

STEP 2. Next we show that the initial function is defined as the $L^2$ limit of $u(t,x)$ as $t \to 0$, that is, we show that $t \to u(\cdot,t) \in C((0,\infty); L^2(D))$ and $u(\cdot,t)$ takes the initial datum $f$ in the sense of $L^2(D)$, that is,

$$\|u(\cdot,t) - f\|_{2,D} \to 0 \quad \text{as } t \to 0.$$

Note $E_\beta(-z)$ is completely monotone for $z \geq 0$ (cf. Agrawal [1]). In particular, $E_\beta(-\lambda_n t^\beta)$ is a nonincreasing function of $\lambda_n$. Hence, we have

$$u(t,x) - f(x) = \sum_{n=1}^{\infty} \bar{f}(n)(E_\beta(-\lambda_n t^\beta) - 1)\phi_n(x).$$

Fix $\epsilon \in (0,1)$ and choose $n_0 = n_0(\epsilon)$ as in (3.8). Then,

$$\|u(\cdot,t) - f\|_{2,D}^2 = \sum_{n=1}^{\infty} (\bar{f}(n))^2 (E_\beta(-\lambda_n t^\beta) - 1)^2$$

$$\leq \sum_{n=1}^{n_0(\epsilon)} (\bar{f}(n))^2 (E_\beta(-\lambda_n t^\beta) - 1)^2$$

$$+ \sum_{n=n_0(\epsilon)+1}^{\infty} (\bar{f}(n))^2 (E_\beta(-\lambda_n t^\beta) - 1)^2$$

$$\leq (1 - E_\beta(-\lambda_{n_0} t^\beta))^2 \|f\|_{2,D} + \epsilon$$

and now the claim follows, since $E_\beta(-\lambda_n t^\beta) \to 1$, as $t \to 0$. The continuity of $t \mapsto u(t, \cdot)$ in $L^2(D)$ at every point $t \in (0, \infty)$ can be proved in a similar fashion.

STEP 3. We now obtain a decay estimate for the solution $u(t, x)$. Using Parseval's identity, the fact that $\lambda_n$ is increasing in $n$, and the fact that $E_\beta(-\lambda_n t^\beta)$ is nonincreasing for $n \geq 1$, we get

$$\|u(\cdot, t)\|_{2,D} \leq E_\beta(-\lambda_1 t^\beta) \|f\|_{2,D}.$$

STEP 4. We show uniform and absolute convergence of the series (3.5) defining $u(t, x)$, which will imply that the classical solution to (3.1) is given by this equation. We do this by showing that (3.6) is a Cauchy sequence in $L^\infty(D)$ uniformly in $t \geq 0$.

Applying the Green's second identity we see that

$$\overline{\Delta f}(n) = \int_D \Delta f(x) \phi_n(x)\, dx = (-\lambda_n) \bar{f}(n).$$

Hence the series defining $\Delta f$ is $\sum_{n=1}^{\infty} -\lambda_n \phi_n(x) \bar{f}(n)$. So, this series is absolutely and uniformly convergent. Let $\epsilon > 0$. Since $\lambda_n \to \infty$ as $n \to \infty$, there exists an $n_0(\epsilon)$ so that for all $x \in D$,

$$(3.10) \qquad \sum_{n=n_0(\epsilon)}^{\infty} |\bar{f}(n)||\phi_n(x)| \leq \sum_{n=n_0(\epsilon)}^{\infty} \lambda_n |\bar{f}(n)||\phi_n(x)| < \epsilon.$$

This is possible since we have assumed that $\Delta f$ has an expansion which is uniformly and absolutely convergent in $L^\infty(D)$. We will freely use the fact that the series defining $f$ also converges absolutely and uniformly.

For $N > M > n_0(\epsilon)$ and $t \geq 0$ and $x \in D$,

$$(3.11) \quad |u_N(t, x) - u_M(t, x)| = \left|\sum_{n=M}^{N} \phi_n(x) \bar{f}(n) E_\beta(-\lambda_n t^\beta)\right|$$

$$\leq \frac{c}{(1 + \lambda_{n_0} t^\beta)} \sum_{n=M}^{M} |\phi_n(x)||\bar{f}(n)| < c\epsilon.$$

This shows that the sequence $u_N(t, x)$ is a Cauchy sequence in $L^\infty(D)$ and so has a limit in $L^\infty(D)$. Hence, the series in (3.5) is absolutely and uniformly convergent. It follows that $u(t, x)$ satisfies the boundary conditions in equation (3.1).

STEP 5. Next we show that the fractional time derivative $\partial_t^\beta$ and the Laplacian $\Delta$ can be applied term by term in (3.5). Note that

$$\left|\sum_{n=1}^{\infty} \bar{f}(n) E_\beta(-\lambda_n t^\beta) \Delta \phi_n(x)\right| = \left|\sum_{n=1}^{\infty} \bar{f}(n) E_\beta(-\lambda_n t^\beta) \lambda_n \phi_n(x)\right|$$

$$\text{(3.12)} \quad \leq \sum_{n=1}^{\infty} |\phi_n(x)||\bar{f}(n)| \frac{\lambda_n}{1+\lambda_n t^\beta}$$

$$\leq t^{-\beta} \sum_{n=1}^{\infty} |\phi_n(x)||\bar{f}(n)| < \infty,$$

where the last inequality follows from the fact that the eigenfunction expansion of $f$ converges absolutely and uniformly. Then the series

$$\sum_{n=1}^{\infty} \bar{f}(n) E_\beta(-\lambda_n t^\beta) \Delta \phi_n(x)$$

is absolutely convergent in $L^\infty(D)$ uniformly in $(t_0, \infty)$ for any $t_0 > 0$. The function $E_\beta(-\lambda t^\beta)$ is an eigenfunction of the Caputo fractional derivative, with $\partial_t^\beta E_\beta(-\lambda t^\beta) = -\lambda E_\beta(-\lambda t^\beta)$ for any $\lambda > 0$, which follows by taking the Laplace transform of both sides of this equation. Hence

$$\sum_{n=1}^{\infty} \bar{f}(n) \phi_n(x) \frac{\partial^\beta}{\partial t^\beta} E_\beta(-\lambda_n t^\beta) = \sum_{n=1}^{\infty} \bar{f}(n) E_\beta(-\lambda_n t^\beta) \Delta \phi_n(x)$$

since the two series are equal term by term, and since the series on the right converges.

Now it is easy to check that the fractional time derivative and Laplacian can be applied term by term in (3.5) to give

$$\frac{\partial^\beta}{\partial t^\beta} u(t,x) - \Delta u(t,x)$$

$$= \sum_{n=1}^{\infty} \bar{f}(n) \left[ \phi_n(x) \frac{\partial^\beta}{\partial t^\beta} E_\beta(-\lambda_n t^\beta) - E_\beta(-\lambda_n t^\beta) \Delta \phi_n(x) \right] = 0,$$

so that the PDE in (3.1) is satisfied. Thus, we conclude that $u$ defined by (3.5) is a classical (strong) solution to (3.1).

Using [21], equation (17),

$$\left| \frac{dE_\beta(-\lambda_n t^\beta)}{dt} \right| \leq c \frac{\lambda_n t^{\beta-1}}{1+\lambda_n t^\beta} \leq c\lambda_n t^{\beta-1},$$

we get

$$\left| \frac{\partial u(x,t)}{\partial t} \right| \leq \sum_{n=1}^{\infty} |\bar{f}(n)| \left| \frac{dE_\beta(-\lambda_n t^\beta)}{dt} \right| |\phi_n(x)|$$

$$\leq ct^{\beta-1} \sum_{n=1}^{\infty} \lambda_n |\bar{f}(n)||\phi_n(x)| = t^{\beta-1} g(x).$$

Since $\Delta f$ has absolutely and uniformly convergent series expansion with respect to $\{\phi_n : n \in \mathbb{N}\}$, we have $g(x) \in L^\infty(D)$.

Now from the results obtained above it follows that $u \in \mathcal{H}_\Delta(D_\infty) \cap C_b(\bar{D}_\infty)$.

We next show that $u \in C^1(\bar{D})$; this follows from the bounds in [16], Theorem 8.33, and the absolute and uniform convergence of the series defining $f$.

(3.13) $$|\phi_n|_{1,\alpha;D} \leq C(1+\lambda_n)\sup_D |\phi_n(x)|,$$

where $C = C(d, \lambda, \Lambda, \partial D)$ and $\lambda$ is the constant in the definition of uniform ellipticity of $L$ in (2.2) and $\Lambda$ is the bound in (2.3). Here

$$|u|_{k,\alpha;D} = \sup_{|\gamma|=k}[D^\gamma u]_{\alpha,D} + \sum_{j=0}^{k}\sup_{|\gamma|=j}\sup_D |D^\gamma u|, \qquad k=0,1,2,\ldots,$$

and

$$[D^\gamma u]_{\alpha,D} = \sup_{x,y\in D, x\neq y} \frac{|D^\gamma u(x) - D^\gamma u(y)|}{|x-y|^\alpha}$$

are norms on $C^{k,\alpha}(\bar{D})$. Hence

$$|u(\cdot,t)|_{1,\alpha;D} \leq C\sum_{n=1}^{\infty} |\bar{f}(n)|E_\beta(-\lambda_n t^\beta)(1+\lambda_n)\sup_D |\phi_n(x)|$$
$$\leq C\sum_{n=1}^{\infty} |\bar{f}(n)|\frac{1+\lambda_n}{1+\lambda_n t^\beta}\sup_D |\phi_n(x)|$$
$$\leq Ct^{-\beta}\sum_{n=1}^{\infty}\sup_D |\phi_n(x)||\bar{f}(n)| + C\sum_{n=1}^{\infty}\sup_D |\phi_n(x)||\bar{f}(n)| < \infty.$$

STEP 6. Next we obtain the stochastic solution to (3.1), by inverting the $\phi_n$-Laplace transform, using a method similar to [24]. Using the fact that for any $\lambda > 0$,

$$\int_0^\infty e^{-\lambda l}\, dl = \frac{1}{\lambda},$$

we get

(3.14)
$$\hat{u}(s,n) = \frac{\bar{f}(n)s^{\beta-1}}{s^\beta + \lambda_n}$$
$$= \bar{f}(n)s^{\beta-1}\int_0^\infty e^{-(s^\beta+\lambda_n)l}\, dl$$
$$= \int_0^\infty e^{-\lambda_n l}\bar{f}(n)s^{\beta-1}e^{-ls^\beta}\, dl.$$

The $\phi_n$-transform of the killed semigroup $T_D(l)f(x) = \sum_{m=1}^{\infty} e^{-\lambda_m l}\phi_m(x) \times \bar{f}(m)$ from (2.5) is found as follows. Since $\{\phi_n, n \in \mathbb{N}\}$ is a complete orthonormal basis of $L^2(D)$, we get

$$\overline{[T_D(l)f]}(n) = \int_D \phi_n(x)T_D(l)f(x)\,dx = \int_D e^{-l\lambda_n}\bar{f}(n)(\phi_n(x))^2\,dx = e^{-l\lambda_n}\bar{f}(n).$$

By [24], equation (8), we have

$$s^{\beta-1}e^{-ls^\beta} = \frac{1}{\beta}\int_0^\infty te^{-st}g_\beta(tl^{-1/\beta})l^{-(1+1/\beta)}\,dt,$$

where $g_\beta(t)$ is the smooth density of the stable subordinator with Laplace transform $\tilde{g}_\beta(s) = \int_0^\infty e^{-st}g_\beta(t)\,dt = e^{-s^\beta}$. Combining the last two results, we get

$$\begin{aligned}
\hat{u}(s,n) &= \int_0^\infty \overline{[T_D(l)f]}(n)\left[\frac{1}{\beta}\int_0^\infty te^{-st}g_\beta(tl^{-1/\beta})l^{-(1+1/\beta)}\,dt\right]dl \\
&= \int_0^\infty e^{-st}\left[\int_0^\infty \overline{[T_D(l)f]}(n)\frac{t}{\beta}g_\beta(tl^{-1/\beta})l^{-(1+1/\beta)}\,dl\right]dt.
\end{aligned} \quad (3.15)$$

Invert the above Laplace transform to obtain

$$\bar{u}(t,n) = \int_0^\infty \overline{[T_D(l)f]}(n)\frac{t}{\beta}g_\beta(tl^{-1/\beta})l^{-(1+1/\beta)}\,dl.$$

Inverting the $\phi_n$-transform is equal to multiplying both sides of the above equation by $\phi_n$ and then summing up from $n = 1$ to $\infty$.

We have uniqueness of this inverse in $L^2(D)$ and $u(t, \cdot) \in L^2(D)$ for each fixed $t \geq 0$. Inverting the $\phi_n$-transform, we obtain

$$u(t,x) = \sum_{n=1}^{\infty} \phi_n(x)\int_0^\infty \overline{[T_D(l)f]}(n)\frac{t}{\beta}g_\beta(tl^{-1/\beta})l^{-(1+1/\beta)}\,dl.$$

Using

$$(3.16) \quad \int_0^\infty e^{-l\lambda_n}\frac{1}{\beta l}tg_\beta(tl^{-1/\beta})l^{-1/\beta}\,dl = \int_0^\infty e^{-\lambda_n(t/l)^\beta}g_\beta(l)\,dl = E_\beta(-\lambda_n t^\beta),$$

we get

$$\begin{aligned}
u(t,x) &= \sum_{n=1}^{\infty}\phi_n(x)\int_0^\infty \overline{[T_D(l)f]}(n)\frac{t}{\beta}g_\beta(tl^{-1/\beta})l^{-(1+1/\beta)}\,dl \\
&= \sum_{n=1}^{\infty}\bar{f}(n)\phi_n(x)E_\beta(-\lambda_n t^\beta).
\end{aligned}$$

After showing the absolute and uniform convergence of the series defining $u$, we can use a Fubini–Tonelli type argument to interchange order of summation and integration in the following to get stochastic representation of the solution as

(3.17)
$$\begin{aligned} u(t,x) &= \sum_{n=1}^{\infty} \phi_n(x) \int_0^{\infty} \overline{[T_D(l)f]}(n) \frac{t}{\beta} g_\beta(tl^{-1/\beta}) l^{-(1+1/\beta)} \, dl \\ &= \int_0^{\infty} \left[ \sum_{n=1}^{\infty} \phi_n(x) \bar{f}(n) e^{-l\lambda_n} \right] \frac{t}{\beta} g_\beta(tl^{-1/\beta}) l^{-(1+1/\beta)} \, dl \\ &= \int_0^{\infty} T_D(l) f(x) \frac{t}{\beta} g_\beta(tl^{-1/\beta}) l^{-(1+1/\beta)} \, dl \\ &= E_x[f(X(E_t)) I(\tau_D(X) > E_t)]. \end{aligned}$$

The last equality follows from a simple conditioning argument, and the density formula (2.7) for the inverse stable subordinator.

STEP 7. Finally we prove uniqueness. Let $u_1, u_2$ be two solutions of (3.1) with initial data $u(0,x) = f(x)$ and Dirichlet boundary condition $u(t,x) = 0$ for $x \in \partial D$. Then $U = u_1 - u_2$ is a solution of (3.1) with zero initial data and zero boundary value. Taking $\phi_n$-transform on both sides of (3.1) we get

$$\frac{\partial^\beta}{\partial t^\beta} \bar{U}(t,n) = -\lambda_n \bar{U}(t,n), \qquad \bar{U}(0,n) = 0,$$

and then $\bar{U}(t,n) = 0$ for all $t > 0$ and all $n \geq 1$. This implies that $U(t,x) = 0$ in the sense of $L^2$ functions, since $\{\phi_n : n \geq 1\}$ forms a complete orthonormal basis for $L^2(D)$. Hence, $U(t,x) = 0$ for all $t > 0$ and almost all $x \in D$. Since $U$ is a continuous function on $D$, we have $U(t,x) = 0$ for all $(t,x) \in [0,\infty) \times D$, thereby proving uniqueness. □

COROLLARY 3.2. *The solution in Theorem (3.1) has also the following representation:*

$$u(t,x) = E_x[f(X(E_t)) I(\tau_D(X) > E_t)] = E_x[f(X(E_t)) I(\tau_D(X(E)) > t)].$$

PROOF. Given a continuous stochastic process $X_t$ on $\mathbb{R}^d$, and an interval $I \subset [0,\infty)$, we denote $X(I) = \{X_u : u \in I\}$. Since the domain $D$ is open and $X_u$ is continuous, it follows that

$$\{\tau_D(X) > t\} = \{X([0,t]) \subset D\}.$$

To see this, first note that $X_u \in D$ for all $0 \leq u \leq t$ implies that $X_t \in D$. Choose an open set $B \subset D$ such that $X_t \in B$. Since $X$ is continuous, the set $X^{-1}(B) = \{u : X_u \in B\}$ is open, and it contains $t$. Then it also contains $t + \varepsilon$ for any $\varepsilon > 0$ sufficiently small, so that $X_u \in D$ for all $0 \leq u \leq t + \varepsilon$, and hence $\tau_D(X) > t$. This shows that $\{X([0,t]) \subset D\} \subseteq \{\tau_D(X) > t\}$. Conversely, if $\tau_D(X) > t$ then $X_u \in D$ for all $0 \le u \le t$, so that $\{\tau_D(X) > t\} \subseteq \{X([0,t]) \subset D\}$, and hence these two sets are equal.

Next note that, since $E_t$ is continuous and monotone nondecreasing, we have $E([0,t]) = [0, E_t]$. To see this, first note that, since $E_t$ is nondecreasing with $E_0 = 0$, $u \in [0,t]$ implies $E_u \in [0, E_t]$, so that $E([0,t]) \subset [0, E_t]$. Conversely, if $y \in [0, E_t]$ then suppose $y \notin E([0,t])$. Then $(E[0,t])$ is not connected, which contradicts the fact that the continuous image of a connected set is connected. Hence $[0, E_t] \subset E([0,t])$, so these two sets are equal.

Finally, we observe that

$$\{\tau_D(X(E)) > t\} = \{X(E([0,t])) \subset D\}$$
$$= \{X([0, E_t]) \subset D\} = \{\tau_D(X) > E_t\}$$

which completes the proof. $\square$

REMARK 3.3. If we subordinate Brownian motion $X$ by a stable subordinator $D_t$ in the sense of Bochner, then the conclusions of Corollary 3.2 do not hold since a stable subordinator does not have continuous sample paths; see, for example, Song and Vondraček [34]. Subordination by the inverse stable subordinator has the property that killing the process $X$ and then subordinating by $E_t$ is the same as subordinating $X$ by $E_t$ and then killing, due to the fact that the sample paths of $E_t$ are continuous and nondecreasing.

COROLLARY 3.4. *Let $f \in C_c^{2k}(D)$ be a $2k$-times continuously differentiable function of compact support in $D$. If $k > 1 + 3d/4$, then the equation (3.1) has a classical (strong) solution. In particular, if $f \in C_c^\infty(D)$, then the solution of equation (3.1) is in $C^\infty(D)$.*

PROOF. By Example 2.1.8 of [12], $|\phi_n(x)| \le (\lambda_n)^{d/4}$. Also, from Corollary 6.2.2 of [13], we have $\lambda_n \sim n^{2/d}$.

Applying the Green's second identity $k$-times, we get

(3.18) $$\overline{\Delta^k f}(n) = \int_D \Delta^k f(x) \phi_n(x)\, dx = (-\lambda_n)^k \bar{f}(n).$$

Using the Cauchy–Schwarz inequality and the fact $f \in C_c^{2k}(D)$, we get

$$\overline{\Delta^k f}(n) \le \left[\int_D (\Delta^k f(x))^2\, dx\right]^{1/2} \left[\int_D (\phi_n(x))^2\, dx\right]^{1/2}$$
$$= \left[\int_D (f^{2k}(x))^2\, dx\right]^{1/2} = c_k,$$

where $c_k$ is a constant independent of $n$.

This and equation (3.18) give $|\bar{f}(n)| \leq c_k(\lambda_n)^{-k}$.
Since

$$\Delta f(x) = \sum_{n=1}^{\infty} -\lambda_n \bar{f}(n) \phi_n(x),$$

to get the absolute and uniform convergence of the series defining $\Delta f$, we consider

$$\sum_{n=1}^{\infty} \lambda_n |\phi_n(x)||\bar{f}(n)| \leq \sum_{n=1}^{\infty} (\lambda_n)^{d/4+1} c_k (\lambda_n)^{-k}$$

$$\leq c_k \sum_{n=1}^{\infty} (n^{2/d})^{d/4+1-k} = c_k \sum_{n=1}^{\infty} n^{1/2+2/d-2k/d}$$

which is finite if $2k/d - 1/2 - 2/d > 1$, that is, $k > 1 + 3d/4$. □

REMARK 3.5. Things are rather easier in an interval $(0, M)$, since eigenfunctions and eigenvalues are explicitly known. Eigenvalues of the Laplacian on $(0, M)$ are $(n\pi/M)^2$, and the corresponding eigenfunctions are $\sin(n\pi x/M)$, for $n = 1, 2, \ldots$. The form of the solution in equation (3.5) on a bounded interval $(0, M)$ in $\mathbb{R}$ was obtained by [1].

Let $\beta \in (0, 1)$, $D_\infty = (0, \infty) \times D$ and define

$$\mathcal{H}_L(D_\infty) \equiv \left\{ u : D_\infty \to \mathbb{R} : \frac{\partial}{\partial t} u(t, x), \frac{\partial^\beta}{\partial t^\beta} u(t, x), Lu(t, x) \in C(D_\infty) \right\},$$

$$\mathcal{H}_L^b(D_\infty) \equiv \mathcal{H}_L(D_\infty)$$
$$\cap \left\{ u : \left|\frac{\partial}{\partial t} u(t, x)\right| \leq g(x) t^{\beta-1}, \text{for some } g \in L^\infty(D), t > 0 \right\}.$$

Next we state and prove an extension of Theorem 3.1 to general uniformly elliptic second-order operators.

THEOREM 3.6. *Let $0 < \alpha < 1$. Let $D$ be a bounded domain with $\partial D \in C^{1,\alpha}$ and $a_{ij} \in C^\alpha(\bar{D})$. Let $\{X_t\}$ be a continuous Markov process corresponding to the generator $L$ defined in (2.1). Let $T_D(t)$ be the killed semigroup corresponding to the process $\{X_t\}$ in $D$. Let $\{E_t\}$ be the process inverse to a stable subordinator of index $\beta \in (0, 1)$ independent of $\{X_t\}$. Let $f \in D(L_D) \cap C^1(\bar{D}) \cap C^2(D)$ for which the eigenfunction expansion (of $Lf$) with respect to the complete orthonormal basis $\{\psi_n : n \geq 1\}$ converges uniformly and absolutely. Then the (classical) solution of*

$$u \in \mathcal{H}_L^b(D_\infty) \cap C_b(\bar{D}_\infty) \cap C^1(\bar{D}),$$

$$\frac{\partial^\beta}{\partial t^\beta} u(t,x) = L u(t,x), \qquad x \in D, \ t \geq 0,$$

(3.19)
$$u(t,x) = 0, \qquad x \in \partial D, \ t \geq 0;$$
$$u(0,x) = f(x), \qquad x \in D,$$

*is given by*

(3.20)
$$\begin{aligned}
u(t,x) &= E_x[f(X(E_t))I(\tau_D(X) > E_t)] \\
&= E_x[f(X(E_t))I(\tau_D(X(E)) > t)] \\
&= \frac{t}{\beta} \int_0^\infty T_D(l) f(x) g_\beta(tl^{-1/\beta}) l^{-1/\beta-1} \, dl \\
&= \int_0^\infty T_D((t/l)^\beta) f(x) g_\beta(l) \, dl.
\end{aligned}$$

PROOF. Suppose $u(t,x) = G(t)F(x)$ is a solution of (3.19). Substituting in the PDE (3.19) leads to

$$F(x) \frac{\partial^\beta}{\partial t^\beta} G(t) = G(t) L F(x)$$

and now dividing both sides by $G(t)F(x)$, we obtain

$$\frac{\partial^\beta}{\partial t^\beta} G(t)/G(t) = \frac{LF(x)}{F(x)} = -\mu.$$

That is,

(3.21) $$\frac{\partial^\beta}{\partial t^\beta} G(t) = -\mu G(t), \qquad t > 0,$$

(3.22) $$LF(x) = -\mu F(x), \qquad x \in D, \ F|_{\partial D} = 0.$$

Problem (3.22) is solved by an infinite sequence of pairs $(\mu_n, \psi_n)$, $n \geq 1$, where $0 < \mu_1 < \mu_2 \leq \mu_3 \leq \cdots$ is a sequence of numbers such that $\mu_n \to \infty$, as $n \to \infty$, and $\psi_n$ is a sequence of functions that form a complete orthonormal set in $L^2(D)$ [cf. (2.4)]. In particular, the initial function $f$ regarded as an element of $L^2(D)$ can be represented as

(3.23) $$f(x) = \sum_{n=1}^\infty \bar{f}(n) \psi_n(x)$$

and by the Parseval identity

(3.24) $$\|f\|_{2,D} = \sum_{n=1}^\infty (\bar{f}(n))^2.$$

Then with $\mu_n$ determined by (3.22) and recalling that $E_\beta(-\mu_n t^\beta)$ solves equation (3.21) with $\mu = \mu_n$, we obtain

$$G(t) = G_0(n)E_\beta(-\mu_n t^\beta),$$

where $G_0(n)$ is selected to satisfy the initial condition $f$. We will show that

(3.25) $$u(t,x) = \sum_{n=1}^\infty \bar{f}(n)E_\beta(-\mu_n t^\beta)\psi_n(x)$$

solves the PDE (3.19).

Define approximate solutions of the form

(3.26) $$u_N(t,x) = \sum_{n=1}^N G_0(n)E_\beta(-\mu_n t^\beta)\psi_n(x), \qquad G_0(n) = \bar{f}(n).$$

STEP 1. The sequence $\{u_N(t,\cdot)\}_{N\in\mathbb{N}}$ is a Cauchy sequence in $L^2(D) \cap L^\infty(D)$, uniformly in $t \in [0,\infty)$. This is proved as in the proof of Theorem 3.1. Hence, the solution to (3.19) is given formally by

(3.27) $$u(t,x) = \sum_{n=1}^\infty \bar{f}(n)E_\beta(-\mu_n t^\beta)\psi_n(x).$$

STEP 2. As in the proof of Theorem 3.1, it can be shown that the series defining $u$ and $Lu$ converge absolutely and uniformly so that we can apply the fractional time derivative and uniformly elliptic operator $L$ term by term to obtain that $u$ defined by (3.27) is a classical solution to (3.19).

STEP 3. The stochastic representation of the solution as given in (3.20) is also valid. We omit the details, as the proof is again similar to that of Theorem 3.1.

STEP 4. We have also a decay estimate for $u$ as in Theorem 3.1, namely,

$$\|u(\cdot,t)\|_{2,D} \le E_\beta(-\lambda_1 t^\beta)\|f\|_{2,D}.$$

Uniqueness is proved as in Theorem 3.1. $\square$

REMARK 3.7. The uniqueness of the solution to the boundary value problem in Theorems 3.1 and 3.6 is shown in [19], Section 4. However, we employ a more direct proof of the uniqueness of the solutions.

REMARK 3.8. It may also be possible to prove the existence of $L^2$-solutions to the fractional Cauchy problem (3.19), as well as the solution formula in the last line of (3.20), using Theorem 3.1 in [5]. To this end, we note that $T_D$, the killed semigroup on the Banach space $L^2(D)$, is uniformly bounded and strongly continuous, since $T_D$ is ultracontractive; see [8, 12]. However, that result does not yield the form of the stochastic solution (3.20) in terms of a killed Markov process, it does not establish classical (strong) solutions, and the proof uses very different methods.

**4. Iterated Brownian motion in bounded domains.** In the previous section, it is shown that a killed Markov process, modified by an inverse $\beta$-stable subordinator, is a stochastic solution to a fractional Cauchy problem on a bounded domain in $\mathbb{R}^d$. The inverse stable subordinators with $\beta = 1/2$ are also related to Brownian subordinators [7]. Since Brownian subordinators are related to higher-order Cauchy problems, this relationship can also be used to connect those higher-order Cauchy problems to their time-fractional analogues. In this section, we establish those connections for Cauchy problems on bounded domains in $\mathbb{R}^d$. We also establish an equivalence between a killed Markov process subordinated to an inverse stable subordinator with $\beta = 1/2$, and the same process subject to a Brownian subordinator. Finally, we identify the boundary conditions that make the two formulations identical. This solves an open problem in [7].

Let $Z_t = X(|Y_t|)$ be the iterated Brownian motion, $D_\infty = (0, \infty) \times D$ and define

$$\mathcal{H}_{\Delta^2}(D_\infty) \equiv \left\{ u : D_\infty \to \mathbb{R} : \frac{\partial}{\partial t} u,\, \Delta^2 u \in C(D_\infty),\, \Delta u \in C^1(\bar{D}),\right.$$
$$\left. \left|\frac{\partial}{\partial t} u(t,x)\right| \leq g(x) t^{-1/2},\, g \in L^\infty(D),\, t > 0 \right\}.$$

Recall that $u \in C^k(\bar{D})$ means that for each fixed $t > 0$, $u(t, \cdot) \in C^k(\bar{D})$.

THEOREM 4.1. *Let $D$ be a domain with $\partial D \in C^{1,\alpha}, 0 < \alpha < 1$. Let $\{X_t\}$ be Brownian motion in $\mathbb{R}^d$, and $\{Y_t\}$ be an independent Brownian motion in $\mathbb{R}$. Let $\{E_t\}$ be the process inverse to a stable subordinator of index $\beta = 1/2$ independent of $\{X_t\}$. Let $f \in D(\Delta_D) \cap C^1(\bar{D}) \cap C^2(D) [\subset L^2(D)]$ be such that the eigenfunction expansion of $\Delta f$ with respect to $\{\phi_n : n \geq 1\}$ converges absolutely and uniformly. Then the (classical) solution of*

$$u \in \mathcal{H}_{\Delta^2}(D_\infty) \cap C_b(\bar{D}_\infty) \cap C^1(\bar{D}),$$

(4.1)
$$\frac{\partial}{\partial t} u(t,x) = \frac{\Delta f(x)}{\sqrt{\pi t}} + \Delta^2 u(t,x), \qquad x \in D,\ t > 0,$$
$$u(t,x) = \Delta u(t,x) = 0, \qquad t \geq 0,\ x \in \partial D,$$
$$u(0,x) = f(x), \qquad x \in D$$

*is given by*

$$u(t,x) = E_x[f(Z_t)I(\tau_D(X) > |Y_t|)] = E_x[f(X(E_t))I(\tau_D(X) > E_t)]$$
(4.2)
$$= E_x[f(X(E_t))I(\tau_D(X(E)) > t)] = 2\int_0^\infty T_D(l)f(x)p(t,l)\,dl,$$

where $T_D(l)$ is the heat semigroup in $D$, and $p(t,l)$ is the transition density of one-dimensional Brownian motion $\{Y_t\}$.

PROOF. Suppose $u$ is a solution to equation (4.1). Taking the $\phi_n$-transform of (4.1) and using Green's second identity, we obtain

$$\text{(4.3)} \qquad \frac{\partial}{\partial t}\bar{u}(t,n) = \frac{-\lambda_n \bar{f}(n)}{\sqrt{\pi t}} + \lambda_n^2 \bar{u}(t,n).$$

Note that the time derivative commutes with the $\phi_n$-transform, as

$$\left|\frac{\partial}{\partial t}u(t,x)\right| \leq g(x)t^{-1/2}, \qquad g \in L^\infty(D),\ t > 0.$$

Taking Laplace transforms on both sides and using the well-known Laplace transform formula

$$\int_0^\infty \frac{t^{-\beta}}{\Gamma(1-\beta)}e^{-st}\,dt = s^{\beta-1}$$

for $\beta < 1$, gives us

$$\text{(4.4)} \qquad s\hat{u}(s,n) - \bar{u}(0,n) = -\lambda_n s^{-1/2}\bar{f}(n) + \lambda_n^2 \hat{u}(s,n).$$

Since $u$ is uniformly continuous on $C([0,\epsilon] \times \bar{D})$, it is also uniformly bounded on $[0,\epsilon] \times \bar{D}$. So, we have $\lim_{t\to 0}\int_D u(t,x)\phi_n(x)\,dx = \bar{f}(n)$. Hence, $\bar{u}(0,n) = \bar{f}(n)$. By collecting the like terms, we obtain

$$\text{(4.5)} \qquad \hat{u}(s,n) = \frac{\bar{f}(n)[1 - \lambda_n s^{-1/2}]}{s - \lambda_n^2}.$$

For fixed $n$ and for large $s$,

$$\text{(4.6)} \quad \hat{u}(s,n) = \frac{\bar{f}(n)[1-\lambda_n s^{-1/2}]}{s-\lambda_n^2} = \frac{s^{-1/2}\bar{f}(n)[s^{1/2}-\lambda_n]}{(s^{1/2}-\lambda_n)(s^{1/2}+\lambda_n)} = \frac{s^{-1/2}\bar{f}(n)}{(s^{1/2}+\lambda_n)}.$$

For any fixed $n \geq 1$, the two formulae (4.5) and (4.6) are well defined and equal for all sufficiently large $s$. Since

$$\text{(4.7)} \qquad \bar{u}(t,n) = \bar{f}(n)E_{1/2}(-\lambda_n t^{1/2}),$$

we can see easily from the proof of Theorem 3.1 that $\bar{u}(t,n)$ is continuous in $t > 0$ for any $n \geq 1$. Hence, the uniqueness theorem for Laplace transforms [4], Theorem 1.7.3, shows that, for each $n \geq 1$, $\bar{u}(t,n)$ is the unique continuous function whose Laplace transform is given by (4.6). Since $x \mapsto u(t,x)$ is an element of $L^2(D)$ for every $t > 0$, and two elements of $L^2(D)$ with the same $\phi_n$-transform are equal $dx$-almost everywhere, we have (4.2) is the unique element of $L^2(D)$ and (4.7) is its $\phi_n$-transform.

Observe that for $\beta = 1/2$, stable subordinator of index $1/2$ has the density [3], Example 1.3.19,

$$g_{1/2}(x) = \frac{1}{\sqrt{4\pi x^3}} \exp\left(-\frac{1}{4x}\right).$$

Therefore the density $p(t,l)$ of $|Y_t|$ is

$$(4.8) \quad 2t g_{1/2}(t/l^2) l^{-3} = \frac{2t}{l^3 \sqrt{4\pi t^3/l^6}} \exp\left(-\frac{l^2}{4t}\right) = \frac{2}{\sqrt{4\pi t}} \exp\left(-\frac{l^2}{4t}\right).$$

As in the proof of Theorem 3.1, we obtain, using (3.16) for $\beta = 1/2$, taking the inverse Laplace transform and then taking the inverse $\phi_n$-transform,

$$(4.9) \quad u(t,x) = \sum_{n=1}^{\infty} \bar{f}(n) E_{1/2}(-\lambda_n t^{1/2}) \phi_n(x)$$

$$= \int_0^{\infty} \left[\sum_{n=1}^{\infty} \bar{f}(n) e^{-\lambda_n l} \phi_n(x)\right] 2t g_{1/2}(t/l^2) l^{-3} \, dl$$

$$(4.10) \quad = \int_0^{\infty} \left[\sum_{n=1}^{\infty} \bar{f}(n) e^{-\lambda_n l} \phi_n(x)\right] p(t,l) \, dl$$

$$= \int_0^{\infty} T_D(l) f(x) p(t,l) \, dl$$

$$(4.11) \quad = E_x[f(Z_t) I(\tau_D(X) > |Y_t|)].$$

Note equation (4.10) follows from equation (4.8), and then equation (4.11) follows by a conditioning argument.

We know from the proof of Theorem 3.1 and Corollary 3.2 that $u$ given in (4.9) can be represented as

$$(4.12) \quad \begin{aligned} u(t,x) &= E_x[f(X(E_t)) I(\tau_D(X) > E_t)] \\ &= E_x[f(X(E_t)) I(\tau_D(X(E)) > t)], \end{aligned}$$

where $E_t$ is the inverse of stable subordinator of index $1/2$.

Now, equations (4.12) and (4.11) lead to

$$u(t,x) = E_x[f(X(E_t)) I(\tau_D(X) > E_t)]$$
$$= E_x[f(X(E_t)) I(\tau_D(X(E)) > t)] = E_x[f(Z_t) I(\tau_D(X) > |Y_t|)].$$

As in the proof of Theorem 3.1, we can easily show that the solution $u(t,x)$ satisfies all the properties in (4.1), except $\Delta^2 u \in C(D_\infty)$ and $\Delta u \in C^1(\bar{D})$. Using (3.13), one can easily show that $\Delta u \in C^1(\bar{D})$. Hence, it suffices to show that $\Delta^2 u \in C(D_\infty)$. To do this, we need only to show the absolute and uniform convergence of the series defining $\Delta^2 u$. To apply $\Delta^2$ term by term to (4.9), we have to show that the series

$$\sum_{n=1}^\infty \bar{f}(n)\phi_n(x)\lambda_n^2 E_{1/2}(-\lambda_n t^{1/2})$$

is absolutely and uniformly convergent for $t > t_0 > 0$.

Note that the $\phi_n$-transform of $\Delta f$ is given by

$$\int_D \phi_n(x)\Delta f(x)\,dx = \int_D f(x)\Delta\phi_n(x)$$
$$= -\lambda_n \int_D f(x)\phi_n(x)\,dx = -\lambda_n \bar{f}(n)$$

and using (3.7), we get

(4.13)
$$\sum_{n=1}^\infty |\bar{f}(n)||\phi_n(x)|\lambda_n^2 E_{1/2}(-\lambda_n t^{1/2}) \le \sum_{n=1}^\infty |\bar{f}(n)||\phi_n(x)|\lambda_n^2 \frac{c}{1+\lambda_n t^{1/2}}$$
$$\le c t_0^{-1/2} \sum_{n=1}^\infty |\bar{f}(n)||\phi_n(x)|\lambda_n < \infty,$$

where the last inequality follows from the absolute and uniform convergence of the eigenfunction expansion of $\Delta f$.

Observe next that the Laplace transform of

$$\frac{\partial}{\partial t}E_\beta(-\lambda_n t^\beta) + \frac{\lambda_n}{\sqrt{\pi t}} - \lambda_n^2 E_\beta(-\lambda_n t^\beta)$$

for $\beta = 1/2$ is

$$\frac{s^{1/2}}{s^{1/2}+\lambda_n} - 1 + \lambda_n s^{-1/2} - \frac{\lambda_n^2 s^{-1/2}}{s^{1/2}+\lambda_n} = 0,$$

since the Laplace transform of $E_{1/2}(-\lambda_n t^{1/2})$ is $\frac{s^{-1/2}}{s^{1/2}+\lambda_n}$ [see equation (3.3)] and the Laplace transform of $1/\sqrt{\pi t}$ is $s^{-1/2}$. Hence, by the uniqueness of Laplace transforms, we get

$$\frac{\partial}{\partial t}E_{1/2}(-\lambda_n t^{1/2}) + \frac{\lambda_n}{\sqrt{\pi t}} - \lambda_n^2 E_{1/2}(-\lambda_n t^{1/2}) = 0.$$

Now applying the time derivative and $\Delta^2$ to the series in (4.9) term by term gives

$$\frac{\partial}{\partial t}u(t,x) - \frac{\Delta f(x)}{\sqrt{\pi t}} - \Delta^2 u(t,x)$$
$$= \sum_{n=1}^{\infty} \bar{f}(n)\left[\phi_n(x)\frac{\partial}{\partial t}E_{1/2}(-\lambda_n t^{1/2}) - \frac{\Delta \phi_n(x)}{\sqrt{\pi t}} - E_{1/2}(-\lambda_n t^{1/2})\Delta^2 \phi_n(x)\right]$$
$$= \sum_{n=1}^{\infty} \bar{f}(n)\phi_n(x)\left[\frac{\partial}{\partial t}E_{1/2}(-\lambda_n t^{1/2}) + \frac{\lambda_n}{\sqrt{\pi t}} - \lambda_n^2 E_{1/2}(-\lambda_n t^{\beta})\right] = 0$$

which shows that the PDE in (4.1) is satisfied. Thus, we conclude that $u$ defined by (4.9) is a classical solution to (4.1).

The uniqueness follows as in the proof of Theorem 3.1. $\square$

REMARK 4.2. Theorem 4.1 also holds with the version of IBM defined by Burdzy [10]. Here, the outer process is the two-sided Brownian motion which is defined by $X_t = X_t^+$ for $t \geq 0$ and $X_t = X_{-t}^-$ for $t \leq 0$, where $X_t^+$ and $X_t^-$ are independent Brownian motions started at 0 and $Y_t$ is an independent one-dimensional Brownian motion. In this case, the IBM takes the form $Z_t = x + X(Y_t)$ and, using a simple conditioning argument, we can show that the function

$$u(t,x) = E_x[f(Z_t)I(-\tau_D(X^-) < Y_t < \tau_D(X^+))]$$

reduces to equation (4.9) and hence is also a solution to (4.1).

REMARK 4.3. Let $f \in C_c^{2k}(D)$ be a $2k$-times continuously differentiable function of compact support in $D$. If $k > 1 + 3d/4$, then the equation (4.1) has a strong solution. In particular, if $f \in C_c^{\infty}(D)$, then the classical solution of equation (4.1) is in $C^{\infty}(D)$. To see this, we use the estimates obtained in Corollary 3.4 to get the absolute and uniform convergence of the series defining $\Delta f$:

$$\sum_{n=1}^{\infty} |\phi_n(x)||\bar{f}(n)|\lambda_n \leq \sum_{n=1}^{\infty} (\lambda_n)^{1+d/4} c(\lambda_n)^{-k} \leq c \sum_{n=1}^{\infty} (n^{2/d})^{1+d/4-k}$$

and this is finite if $k > 1 + 3d/4$.

Recall $D_\infty = (0,\infty) \times D$ and define

$$\mathcal{H}_{L^2}(D_\infty) \equiv \left\{u : D_\infty \to \mathbb{R} : \frac{\partial}{\partial t}u,\ L^2 u \in C(D_\infty),\ Lu \in C^1(\bar{D}),\right.$$
$$\left.\left|\frac{\partial}{\partial t}u(t,x)\right| \leq g(x)t^{-1/2}, g \in L^\infty(D), t > 0\right\}.$$

THEOREM 4.4. *Let $D \subset \mathbb{R}^d$ be a domain with $\partial D \in C^{1,\alpha}, 0 < \alpha < 1$, and $a_{ij} \in C^\alpha(\bar{D})$. Let $\{X_t\}$ be the process in $\mathbb{R}^d$ corresponding to $L$ defined in (2.1), and $\{Y_t\}$ be an independent Brownian motion in $\mathbb{R}$. Let $\{E_t\}$ be the process inverse to a stable subordinator of index $\beta = 1/2$ independent of $\{X_t\}$. Let $f \in D(L) \cap C^1(\bar{D}) \cap C^2(D)[\subset L^2(D)]$ be such that the eigenfunction expansion of $Lf$ with respect to $\{\psi_n : n \geq 1\}$ converges absolutely and uniformly. Then, the (classical) solution of*

(4.14)
$$u \in \mathcal{H}_{L^2}(D_\infty) \cap C_b(\bar{D}_\infty) \cap C^1(\bar{D}),$$
$$\frac{\partial}{\partial t} u(t,x) = \frac{Lf(x)}{\sqrt{\pi t}} + L^2 u(t,x), \qquad x \in D, t > 0,$$
$$u(t,x) = Lu(t,x) = 0, \qquad t \geq 0,\ x \in \partial D,$$
$$u(0,x) = f(x), \qquad x \in D$$

*is given by*

$$u(t,x) = E_x[f(X(|Y_t|))I(I(\tau_D(X) > |Y_t|))] = E_x[f(X(E_t))I(\tau_D(X) > E_t)]$$
$$= E_x[f(X(E_t))I(\tau_D(X(E)) > t)] = 2 \int_0^\infty T_D(l) f(x) p(t,l)\, dl,$$

*where $T_D(l)$ is the killed semigroup in $D$ corresponding to $L$ and $p(t,l)$ is the transition density of one-dimensional Brownian motion $\{Y_t\}$.*

PROOF. We can show that the solution $u$ is given by the series

$$\sum_{n=1}^\infty \bar{f}(n) E_{1/2}(-\mu_n t^{1/2}) \psi_n(x).$$

Similar to the proofs of Theorems 3.1 and 4.1, we can prove that this series can be differentiated term by term to give classical solutions.

Let $u_1, u_2$ be two solutions with initial value $f$ and boundary values $0$. Then $u = u_1 - u_2$ is a solution with initial value equal to $0$.

Use the fact that $L$ is self-adjoint [same as integration by parts since $Lu, u \in C^1(\bar{D})$] to get

$$\int_D \psi_n(x) Lu(t,x)\, dx = \int_D u(t,x) L\psi_n(x)\, dx$$
$$= -\mu_n \int_D u(t,x) \psi_n(x)\, dx = -\mu_n \bar{u}(t,n).$$

Using $Lu(t,x) = 0$ and applying integration by parts twice, we obtain

$$\frac{\partial}{\partial t} \bar{u}(t,n) = \mu_n^2 \bar{u}(t,n).$$

Since the initial function is zero, we obtain that $\bar{u}(t,n) = 0$ for all $t > 0$ and for all $n \geq 1$. This implies that $u(t,x) = 0$ in the sense of $L^2$ functions due to the fact that $\{\psi_n : n \geq 1\}$ forms a complete orthonormal basis of $L^2(D)$. Hence, $u(t,x) = 0$ for all $t > 0$ and almost all $x \in D$. Since $u$ is a continuous function on $D$, we have $u(t,x) = 0$ for all $(t,x) \in [0,\infty) \times D$, thereby proving the uniqueness. $\square$

COROLLARY 4.5. *The two types of PDEs in (3.19) and (4.14), with the respective initial and boundary conditions, are equivalent.*

REMARK 4.6. We discovered in [7] that the two types of PDEs (3.19) for $\beta = 1/2$ and (4.14) are not equivalent without any further assumptions on the boundary behavior of the solutions in strict subdomains of $\mathbb{R}^d$. We get the equivalence of the PDEs (3.19) and (4.14) in bounded domains, if we consider the right conditions on boundary behavior of solutions.

REMARK 4.7. It would also be interesting to consider fractional Cauchy problems on bounded domains with $1 < \beta < 2$, and $\partial_t u(t,x) = 0$ on $x \in \partial D$. For unbounded domains, some results along these lines can be found in [6, 9].

M. M. MEERSCHAERT
DEPARTMENT OF STATISTICS
  AND PROBABILITY
MICHIGAN STATE UNIVERSITY
EAST LANSING, MICHIGAN 48823
USA
E-MAIL: mcubed@stt.msu.edu
URL: http://www.stt.msu.edu/~mcubed/

E. NANE
DEPARTMENT OF MATHEMATICS
  AND STATISTICS
AUBURN UNIVERSITY
221 PARKER HALL
AUBURN, ALABAMA 36849
USA
E-MAIL: nane@auburn.edu

P. VELLAISAMY
DEPARTMENT OF MATHEMATICS
INDIAN INSTITUTE OF TECHNOLOGY
  BOMBAY
POWAI, MUMBAI 400076
INDIA
E-MAIL: pv@math.iit.ac.in